\documentclass{amsart}
\usepackage{amssymb,amsmath,latexsym,amsfonts,amsthm,mathrsfs}
\usepackage{enumerate}
\newtheorem{thm}{Theorem}

\newtheorem*{claim}{Claim}
\newtheorem{claimm}{Claim}

\newtheorem{lemma}{Lemma}

\newtheorem{cor}{Corollary}

\newcommand{\R}{\mathbb{R}}
\newcommand{\Q}{\mathbb{Q}}
\newcommand{\funct}[2]{#1 \longrightarrow #2}
\newcommand{\ot}[1]{\textbf{#1},<^{\textbf{#1}}_{lex}}
\newcommand{\om}[1]{\textbf{#1},<^{\textbf{#1}}}
\newcommand{\m}[1]{\textbf{#1}}
\newcommand{\oc}[1]{\widetilde{\textbf{#1}},<^{\widetilde{\textbf{#1}}}}
\newcommand{\mc}[1]{\widetilde{\textbf{#1}}}

\newcommand{\otc}[1]{\widetilde{\textbf{#1}},<^{\widetilde{\textbf{#1}}}_{lex}}
\newcommand{\restrict}[2]{#1 \upharpoonright #2}
\newcommand{\arrows}[3]{\longrightarrow {#1}^{#2}_{#3}}
\newcommand{\U}{\mathcal{U}_S}
\newcommand{\UU}{\mathcal{U}^{c<}_S}

\newcommand{\Ur}{\textbf{Q}_{S}}

\newcommand{\cUr}{\textbf{U}_{S}}
\newcommand{\iso}{\mathrm{iso}}
\newcommand{\cLO}{\mathrm{cLO}}

\author{Lionel Nguyen Van The}

\address{Equipe de Logique Math\'ematique, UFR de Math\'ematiques,
 (case 7012), Universit\'e Denis Diderot Paris 7, 2 Place Jussieu,
 75251 Paris Cedex 05, France.}

\email{nguyenl@logique.jussieu.fr}

\title{Ramsey degrees of finite ultrametric spaces, ultrametric Urysohn spaces and dynamics of their isometry groups}

\date{June 17, 2004}

\begin{document}

\maketitle

\begin{abstract}
We study Ramsey-theoretic properties of several natural classes of
finite ultrametric spaces, describe the corresponding Urysohn spaces
and compute a dynamical invariant attached to their isometry groups. 
\end{abstract}

\section{Introduction}

The purpose of this note is the study of Ramsey-theoretic properties of various
classes $\mathcal{K}$ of finite ultrametric spaces together with some
dynamical consequences on certain isometry groups. The original motivation of
this work comes from the recent work of Kechris, Pestov and Todorcevic
\cite{KPT} where a general problem is posed. Namely, what are the
classes of finite ordered metric spaces (or more generally the classes
of finite ordered structures) which satisfy both Ramsey and Ordering
properties ? To see what this problem asks in a specific context,
recall that a metric space
$\m{X} = (X, d^\m{X})$ is \textit{ultrametric} when given any $x,
y, z$ in $\m{X}$, \begin{center}
$d^\m{X}(x,z) \leqslant \max(d^\m{X}(x,y), d^\m{X}(y,z))$ \end{center} Given ultrametric
spaces $\m{X}$, $\m{Y}$ and $\m{Z}$, we write $\m{X} \cong \m{Y}$ when there is an isometry from $\m{X}$ onto
$\m{Y}$ and define the set $\binom{\m{Z}}{\m{X}}$ as \begin{center}
 $ \binom{\m{Z}}{\m{X}} = \{ \mc{X} \subset \m{Z} : \mc{X} \cong \m{X}
  \}$ \end{center} For $k,l \in \omega \smallsetminus \{ 0 \}$ and a
triple $\m{X}, \m{Y}, \m{Z}$ of ultrametric spaces, the
symbol \begin{center} $\m{Z} \arrows{(\m{Y})}{\m{X}}{k,l}$ \end{center}
is an abbreviation for the statement: \begin{center} For any $\chi :
\funct{\binom{\m{Z}}{\m{X}}}{k}$
there is $\widetilde{\m{Y}} \in \binom{\m{Z}}{\m{Y}}$ such
that $\chi$ does not take more than $l$ values on
$\binom{\widetilde{\m{Y}}}{\m{X}}$. \end{center} When $l = 1$, this is simply
written $\m{Z} \arrows{(\m{Y})}{\m{X}}{k}$. Given a class $\mathcal{K}$ of
ultrametric spaces and $\m{X} \in \mathcal{K} $, if there is $l
\in \omega \smallsetminus \{ 0 \}$ such that \begin{center} For any $\m{Y}
\in \mathcal{K}$, and any $k \in \omega \smallsetminus \{ 0 \}$,
there exists $\m{Z} \in \mathcal{K} $ such that $\m{Z}
\arrows{(\m{Y})}{\m{X}}{k,l}$ \end{center} we write
$\mathrm{t}_{\mathcal{K}}(\m{X})$ for the least such number
$l$. $\mathrm{t}_{\mathcal{K}}(\m{X})$ is called the \textit{Ramsey degree of} $\m{X}$ \textit{in} $\mathcal{K}$ (this is part of the more general notion of Ramsey
degree for an arbitrary class of structures that has already been
studied in the literature, see for example \cite{F} or \cite{KPT}). It
turns out that a positive answer to the question of \cite{KPT} is
equivalent to the existence and the computation of the Ramsey degree
for every member of a given class of metric spaces (resp. of finite
structures). In this note, we are able to do this for the class of
\textit{finite convexly ordered ultrametric spaces}.

For an ultrametric space $\m{X}$, let $\mathrm{iso}(\m{X})$
denote the set of all isometries from $\m{X}$ into itself, and
$\mathrm{cLO}(\m{X})$ the set of all convex linear orderings
of $X$ (a linear ordering $<$ on $\m{X}$ is \textit{convex} when all the metric balls of
$\m{X}$ are $<$-convex). For $S \subset ]0, + \infty [$, let $\U$ denote the class of all finite
ultrametric spaces with distances in $S$. Let also $\UU$ denote the class of
all finite convexly ordered ultrametric spaces with distances in $S$. 

\begin{thm}
Let $S \subset ]0, + \infty [$. Then every element $\m{X}$ of $\U$ has a
Ramsey degree in $\U$ which is equal to $|\mathrm{cLO}(\m{X})|/|\mathrm{iso}(\m{X})|$.
\end{thm}

It turns out that this sort of results is closely related to purely Ramsey-theoretic
results for some classes of \textit{ordered ultrametric spaces},
that is structures of the form $(\om{Z}) = (Z,
d^\m{Z}, <^\m{Z} )$ where $\m{Z}$ is an ultrametric space and $<^\m{Z}$ is a linear ordering on
$Z$. For two ordered ultrametric spaces $(\om{X})$
and $(\om{Y})$, $(\om{X}) \cong (\om{Y})$ means that there is an
order preserving isometry from $(\om{X})$ to
$(\om{Y})$. The notions \begin{center} $\binom{\om{Z}}{\om{X}}$, $(\om{Z})
\arrows{(\om{Y})}{(\om{X})}{k,l}$ and $(\om{Z})
\arrows{(\om{Y})}{(\om{X})}{k}$ \end{center} are defined along the same lines as in
the unordered case. Now, given a class $\mathcal{K}$ of finite ordered
ultrametric spaces, say that $\mathcal{K}$ has the \textit{Ramsey
  property} when for every
$(\om{X}), (\om{Y}) \in \UU$
and every $k \in \omega \smallsetminus \{ 0 \}$, there is
$(\om{Z}) \in \UU$ such that \begin{center} $ (\om{Z})
\arrows{(\om{Y})}{(\om{X})}{k}$ \end{center} On the other hand, say
that $\mathcal{K}$ has the \textit{ordering property} when given
$\m{X} \in \U$, there is $\m{Y} \in \U$ such that given any convex
orderings $<^{\m{X}}$ and $<^{\m{Y}}$ on $\m{X}$ and $\m{Y}$ respectively, $(\om{Y})$ contains an isomorphic
copy of $(\om{X})$. 

Theorem 1 is linked to the
following results:

\begin{thm}

\label{thm:2}

Let $S \subset ]0, + \infty [$. Then $\UU$ has the Ramsey property.
\end{thm}

\begin{thm}

\label{thm:3}

Let $S \subset ]0, + \infty [$. Then $\UU$ has the ordering property.

\end{thm}

Together with the ones obtained by Ne\v{s}et\v{r}il in \cite{N1},
\cite{N2} (who was also
motivated by the general problem posed in \cite{KPT}) in the context of
finite ordered metric spaces,
these results provide some of the very few known examples of Ramsey classes
of finite ordered metric spaces.  

Given a countable class $\mathcal{K}$ of finite metric spaces, one may
also ask for the existence of a corresponding \textit{Urysohn
space}, that is a countable metric space whose family of finite
subspaces is exactly $\mathcal{K}$ and which is \textit{ultrahomogeneous}
i.e. where any
isometry between finite subspaces can be extended to an isometry of
the whole space (for a survey, see \cite{B2}). For $\mathcal{K} =
\U$, there is, up to isometry, a unique such object (denoted here
$\Ur$) which is well-known (see \cite{B1}, \cite{D}, \cite{P} or more recently
\cite{GK}). Particularly convenient frameworks for this study are the
ones developed by Fra\"iss\'e in the general case (for references on Fra\"iss\'e
theory, see \cite{Fr} or \cite{H}) and, before him, Urysohn (\cite{U})
in the case of metric spaces. Here, we present a simple description of $\Ur$ as well as
its completion $\cUr$ and provide an explicit computation of the universal minimal
flow of the corresponding isometry groups. Recall that for a topological group $G$, a \textit{compact
minimal $G$-flow} is a compact Hausdorff space $X$ together with a
continuous action of $G$ on $X$ for which the orbit of every point is
dense in $X$. It is a general result in topological dynamics that
every topological group $G$ has a compact minimal $G$-flow $M(G)$
which is, moreover, \textit{universal} in the
sense that it can be mapped homomorphically onto any other compact
minimal $G$-flow. Here, we prove:

\begin{thm}    
Let $S \subset ]0, + \infty [$ be countable. Then $M(\iso(\Ur)) = \cLO(\Ur)$
together with the natural action $\funct{\iso (\Ur) \times \cLO (\Ur)}{\cLO (\Ur)}$, $(g,<)
\longmapsto <^g$ defined by $x <^g y$ iff $g^{-1}(x) < g^{-1}(y)$.
Similarly,  $M(\iso(\cUr)) = \cLO(\cUr)$ together with the action $\funct{\iso (\cUr) \times \cLO (\cUr)}{\cLO (\cUr)}$, $(g,<)
\longmapsto <^g$ defined by $x <^g y$ iff $g^{-1}(x) < g^{-1}(y)$.
\end{thm}

This is obtained thanks to the technique used in \cite{KPT} to show
that for the rational Urysohn space $\textbf{U} _0$ corresponding to the class of all finite
metric spaces with rational distances, the universal minimal
flow is the space of all linear orderings on $\textbf{U} _0$.

Finally, before diving into the heart of the matter, there are three
persons I would like to thank sincerely. Alexander Kechris for the
helpful comments and references concerning section 5 and 6. Jordi Lopez Abad for
the numerous, long and colored discussions. And finally, Stevo Todorcevic for
his clear guidance and unalterable patience. 

\section{Trees and finite convexly ordered ultrametric spaces}
In this section, we present a duality between trees and ultrametric
spaces. This fact is the key for most of the proofs in this note. 

A \emph{tree} $\m{T} = (T,<^{\m{T}}) $ is a partially
ordered set such that given any element $t
\in T$, the set $\{ s \in T : s <^{\m{T}} t \}$ is
$<^{\m{T}}$-well-ordered. When every element of $T$ has finitely many
$<^{\m{T}}$-predecessors, $\mathrm{ht}(t) = |\{ s \in T : s <^{\m{T}}
t \}|$ is the \emph{height of} $t
\in \m{T}$ and when $n< \mathrm{ht}(t)$, $t(n)$ denotes the unique
predecessor of $t$ with height $n$. The $m$-th level of $\m{T}$
is $\m{T}(m) = \{t \in T : \mathrm{ht}(t) = m \}$ and the \emph{height
of} $\m{T}$ is the least $m$ such that $\m{T}(m) = \emptyset$. When $|\m{T}(0)| = 1$, we say that $\m{T}$ is \textit{rooted}. When $\m{T}$ is rooted and $s, t \in \m{T}$,
$\Delta (s,t)$ is the largest $n$ such that $s(n) = t(n)$. A linear
ordering $<$ on $\m{T}$ is \emph{lexicographical} when $\m{T}(0) <
\m{T}(1) < \ldots < \m{T}({\mathrm{ht}(\m{T})})$ and given $s, t \in
\m{T}$ in a same level, $s < t$ whenever $s(\Delta(s,t)+1) < t(\Delta(s,t)+1)$. 
From now on, all the trees $\m{T}$ we will consider will be of
finite height, rooted, lexicographically ordered by some ordering $<^{\m{T}}_{lex}$, and
the set $\m{T}^{max}$ of 
all $<^{\m{T}}$-maximal elements of
$\m{T}$ will coincide with the top level set of $\m{T}$.
Given such a tree of height $n$ and a finite sequence
$a_0 > a_1> \ldots >a_{n-1}$ of strictly positive real numbers, there is a
natural convexly ordered ultrametric space structure on
$\m{T}^{max}$ if the distance $d$ is defined by $d(s,t) =
a_{\Delta(s,t)}$. Conversely, given any
convexly ordered ultrametric space $(\om{X})$ with finitely many distances
given by $a_0 > a_1> \ldots >a_{n-1}$, there is a tree $(\ot{T})$ of
height $n$ such that $(\om{X})$ is the natural convexly ordered
ultrametric space associated to $(\ot{T})$ and $(a_i)_{i<n}$. The
elements of $\m{T}$ are the ordered pairs of the form $\left\langle m, b \right\rangle$ where $m \in n$ and  
$b$ is a ball of the form $\{ y \in \m{X} :
d^{\m{X}}(y,x) \leqslant a_m \}$ with $x \in \m{X}$. The structural ordering $<^{\m{T}}$ is
given by $\left\langle l, b \right\rangle <^{\m{T}} \left\langle m, c \right\rangle$ iff ($l < m $ and $ b \subset c$) and the
lexicographical ordering is defined levelwise by $\left\langle m, b \right\rangle
<^{\m{T}}_{lex} \left\langle m, c\right\rangle$ iff  there is $x \in b$ and $y \in c$ such
that $x<^{\m{X}}y$.

\section{Ramsey degrees for finite ultrametric spaces}

In this section, $S \subset ]0, + \infty [$. For $\m{X} \in
\U$, let $\tau (\m{X}) := |\mathrm{cLO}(\m{X})|/|\mathrm{iso}(\m{X})|$. $\tau
(\m{X})$ is essentially the number of all
nonisomorphic structures one can get by adding a convex linear ordering on
$\m{X}$. Indeed, if $<_1, <_2$ are convex linear orderings on
$\m{X}$, then $(\m{X},<_1)$ and $(\m{X},<_2)$ are
isomorphic as finite ordered metric spaces if and only if the
unique order preserving bijection from $(\m{X},<_1)$ to
$(\m{X},<_2)$ is an isometry. This defines equivalence
relation on the set of all finite convexly ordered ultrametric spaces obtained
by adding a convex linear ordering on $\m{X}$. In what follows, an
\textit{order type for} $\m{X}$ is an equivalence class
corresponding to this relation. In this section, we prove
theorem 1, that is we show that any $\m{X} \in \U$ has a
Ramsey degree in $\U$ which is equal to $\tau (\m{X})$. This
result will be obtained via theorem 2 (for which the proof is postponed to
section 4) and theorem 3.

\subsection{The existence result - Getting an upper bound}

\begin{thm}
Let $\m{X} \in \U$. There is $l \in \omega$ such that for
every $\m{Y} \in \U $, $k \in \omega \smallsetminus \{ 0 \}$,
there is $\m{Z} \in \U$ such that $\m{Z}
\arrows{(\m{Y})}{\m{X}}{k,l}$.
\end{thm}

In fact, we prove that the value $l = \tau (\m{X})$ works.
This will provide the result, as well as an upper bound for
$\mathrm{t}_{\U}(\m{X})$.

\begin{proof}
Let $\{ <_{\alpha} : \alpha \in A \}$ be a set of convex linear orderings
on $\m{X}$ such that for every convex linear ordering $<$ on
\m{X}, there is a unique $\alpha \in A$ such that
$(\m{X},<)$ and $(\m{X},<_{\alpha})$ are isomorphic as
finite ordered metric spaces. Then $A$ has size
$\tau (\m{X})$ so without loss of generality, $A = \{1,
\ldots, \tau (\m{X}) \}$. Now, let $<^\m{Y}$ be any
convex linear ordering on $Y$. By Ramsey property
for $\UU$ (theorem \ref{thm:2}) we can find $(\m{Z}_1,<^{\m{Z}_1} ) \in \UU$ such that
$(\m{Z}_1,<^{\m{Z}_1})
\arrows{(\m{Y},<^\m{Y})}{(\m{X},<_1)}{k} $. Now,
construct inductively $(\m{Z}_2,<^{\m{Z}_2}), \ldots,
(\m{Z}_{\tau (\m{X})},<^{\m{Z}_{\tau (\m{X})}})
\in \UU$ such that for every $n \in \{1, \ldots,
\tau (\m{X})-1 \}$, $(\m{Z}_{n+1},<^{\m{Z}_{n+1}})
\arrows{(\m{Z}_n,<^{\m{Z}_n})}{(\m{X},<_{n+1})}{k}$.
Finally, let $\m{Z}=\m{Z}_{\tau (\m{X})}$. Then one can check
that $\m{Z} \arrows{(\m{Y})}{\m{X}}{k,\tau (\m{X})}$.

\end{proof}

\subsection{Ordering property - Reaching the Ramsey degree}

In the previous subsection, reducing the number of
values of a given coloring is possible because the Ramsey
property for $\UU$ allows to color copies of $\m{X}$
according to their order type in $\m{Y}$. Consequently, the
fact that this reduction is not possible under a certain point
should mean that given some $\m{X} \in \U$, it is possible to
find $\m{Y} \in \U$ such that given any convex linear ordering $<$
on $\m{Y}$, every order type of $\m{X}$ is represented
in $(\m{Y},<)$. In this section, we show that this is indeed
the case, which proves theorem \ref{thm:3}. We begin with a simple observation coming from the tree
representation of elements of $\UU$. The proof is left to the reader. 

\begin{lemma}[Reasonability property for $\UU$]

\label{lemma:1}

Given $\m{X} \subset \m{Y}$ in $\U$ and $<^{\m{X}}$ a convex
linear ordering on $\m{X}$, there is a convex linear ordering $<^{\m{Y}}$ on $\m{Y}$
such that $\restrict{ <^{\m{Y}} }{\m{X} } = <^{\m{X}}$.
\end{lemma}

Call an element $\m{Y}$ of $\U$ \textit{convexly order-invariant} when
$(\m{Y},<_1) \cong (\m{Y},<_2)$ whenever $<_1 , <_2$ are convex
linear orderings on $\m{Y}$.
The following result is a direct consequence of the previous lemma: 

\begin{lemma}
Let $\m{X} \subset \m{Y}$ in $\U$, and assume that $\m{Y}$ is convexly 
order-invariant. Then given any convex linear ordering $<$ on $\m{Y}$,
every order type of $\m{X}$ is represented in $(\m{Y}, <)$.
\end{lemma}

\begin{proof}
Let $<$ and $<^{\m{X}}$, be convex linear orderings on $\m{Y}$ and $\m{X}$
respectively. 
Let $<^{\m{Y}}$ be as in the previous lemma. Then $(\om{X})$ is
represented in $(\om{Y}) \cong (\m{Y}, <)$.
\end{proof}

We now show that any element of $\U$ embeds into a convexly order-invariant
one.

\begin{lemma}
Let $\m{X} \in \U$. Then $\m{X}$ embeds into a convexly order-invariant $\m{Y} \in \U$. 
\end{lemma}

\begin{proof}
Let $a_0 > a_1 > \ldots > a_{n-1}$ enumerate the distances appearing
in $\m{X}$. The tree representation of $\m{X}$ has $n$ levels. Now,
observe that such a tree can be embedded into a tree of height $n$ and
where all the nodes of a same level have the same number of
immediate successors, and that the ultrametric space associated to that
tree is convexly order-invariant. 
\end{proof}

Combining these results, we get: 

\begin{cor}

\label{cor:1}

Let $X \in \U$. There is $Y \in \U$ such that given any convex linear
ordering $<^{\m{Y}}$ on $Y$, every order type of $X$ is represented in $(\om{Y})$.
\end{cor}

Theorem \ref{thm:3} follows then directly. We can now connect $\tau(\m{X})$ and
$\mathrm{t} _{\U} (\m{X})$.

\begin{thm}
Let $\m{X} \in \U$. Then there is $\m{Y} \in \U$ such
that for every $\m{Z} \in \U$, there is $\chi :
\funct{\binom{\m{Z}}{\m{X}}}{\tau (\m{X})}$ with the
property: Given any $\widetilde{\m{Y}} \in
\binom{\m{Z}}{\m{Y}} $, $\chi$ takes $\tau (\m{X})$
values on $\binom{\widetilde{\m{Y}}}{\m{X}}$.
\end{thm}

\begin{proof}
Fix $\m{X} \in \U$ and let $\m{Y} \in \U$ be as in corollary
\ref{cor:1}: For any convex linear ordering $<$ on $\m{Y}$,
$(\m{Y},<)$ contains a copy of each order type of
$\m{X}$. Now, let $\m{Z} \in \U$ and pick $<^\m{Z}$
any convex linear ordering on $\m{Z}$. Define a coloring $\chi :
\funct{\binom{\m{Z}}{\m{X}}}{\tau (\m{X})}$ which colors any copy
$\widetilde{\m{X}}$ of $\m{X}$ according to the order type
of $(\widetilde{\m{X}},\restrict{<^\m{Z}}{\widetilde{\m{X}}})$. Now, if
possible, let $\widetilde{\m{Y}} \in \binom{\m{Z}}{\m{Y}} $.
Then $(\widetilde{\m{Y}}, \restrict{<^\m{Z}}{\widetilde{\m{Y}}}) $
contains a copy of every order type of $\m{X}$, and

\[
|\chi '' \binom{\widetilde{\m{Y}}}{\m{X}} | = \tau (\m{X})
\]

\end{proof}

A direct consequence of this theorem is:

\begin{cor}

\label{cor:6}

For any $\m{X} \in \U$, $\mathrm{t}_{\U} (\m{X}) =
\tau(\m{X})$.
\end{cor}

At that point, some comments are of interest. The first one concerns a
parallel with the results related to the class $\mathcal{M}$ of finite metric spaces. 
Indeed, Ne\v{s}et\v{r}il proved that the class $\mathcal{M} ^<$ of all
finite ordered metric spaces has the Ramsey
property (\cite{N1}) as well as the ordering
property (\cite{N2}). Hence, every element $\m{X}$ in $\mathcal{M}$
has a Ramsey degree which is equal to
$|\mathrm{LO}(\m{X})|/|\mathrm{iso}(\m{X})|$ ($\mathrm{LO}(\m{X})$ being the
set of all linear orderings on $\m{X}$). This fact has two
consequences. On the one hand, the only Ramsey objects (i.e. those for which
$\mathrm{t}_{\mathcal{M}}(\m{X}) = 1$) are the equilateral ones. On the other
hand, there are objects for which the Ramsey degree is $\mathrm{LO}(\m{X})$ (i.e. $|\m{X}|!$), those for which there
is no nontrivial isometry. In the present case, the situation is a bit different: First, 
the ultrametric spaces for which the true Ramsey property holds are those for which the corresponding
tree is uniformly branching on each level. Hence, in the class $\U$,
every element can be embedded into a Ramsey object, a fact which does
not hold in the class of all finite metric spaces. Second, one can
notice that any finite ultrametric space has a nontrivial isometry
(this fact is obvious via the tree representation). Thus, the Ramsey
degree of $\m{X}$ is always strictly less than
$|\mathrm{cLO}(\m{X})|$. In fact, a simple computation shows that the
highest value $\mathrm{t}_{\U}(\m{X})$ can get if the size of
$\m{X}$ is fixed is $2^{|\m{X}|-2}$ and is realized when the tree
associated to $\m{X}$ is a comb, i.e. all the branching nodes are
placed on a same branch. 

The second comment concerns the class $\mathcal{U} ^< _S$ of all
finite ordered ultrametric spaces. One can show that for this class,
Ramsey property implies Ordering property. As a consequence,
$\mathcal{U} ^< _S$ cannot have the Ramsey property. Here is a simple
combinatorial argument for that: Otherwise, any
element in $\U$ would have a Ramsey degree equal to
$|\mathrm{LO}(\m{X})|/|\mathrm{iso}(\m{X})|$, a contradiction with corollary \ref{cor:6}. However, let us also 
mention that a much deeper reason is hidden behind that fact (see
\cite{KPT}, section 9).

\section{Ramsey property for finite convexly ordered ultrametric spaces}

The purpose of this section is to provide a proof of theorem 2.
To do that, let us introduce some notations for the partition
calculus on trees.
Given trees $(\ot{T})$ and $(\ot{S})$ as described in section 2, say that
they are \textit{isomorphic} when there is a bijection between them which
preserves both the structural and the lexicographical orderings.
Also, given a tree $(\ot{U})$, write $\binom{\ot{U}}{\ot{T}}$ for the
set $\{ (\otc{T}) : \widetilde{\m{T}} \subset \m{U} \wedge (\otc{T}) \cong (\ot{T}) \}$.
Now, if $(\ot{S}), (\ot{T})$ and $(\ot{U})$ are trees, let $(\ot{U})
\arrows{(\ot{T})}{(\ot{S})}{k}$ abbreviate the statement:  
For any $\chi : \funct{\binom{\ot{U}}{\ot{S}}}{k}$ there is
$(\otc{T}) \in \binom{\ot{U}}{\ot{T}}$ such that $\chi$ is
constant on $\binom{\otc{T}}{\ot{S}}$.

\begin{lemma}

\label{lemma:2}

Given an integer $k \in \omega \smallsetminus \{ 0 \}$, a finite tree $(\ot{T})$ and a subtree $(\ot{S})$ of $(\ot{T})$ such that
$\mathrm{ht}(\m{T}) = \mathrm{ht}(\m{S})$, there is a finite tree $(\ot{U})$ such that $\mathrm{ht}(\m{U}) = \mathrm{ht}(\m{T})$
and $(\ot{U}) \arrows{(\ot{T})}{(\ot{S})}{k}$.
\end{lemma}

A natural way to proceed is by induction on
$\mathrm{ht}(\m{T})$. Actually, it is so natural that after having done so, we realized that this method had already been
used in \cite{F} where the exact same
result is obtained ! Consequently, we choose to provide a different proof
which uses the notion of ultrafilter-tree.

\begin{proof}

For the sake of clarity, we sometimes not mention the lexicographical
orderings. For example, $\m{T}$ stands for $(\ot{\m{T}})$. So let $\m{T} \subset \m{S}$ be some finite trees of height $n$ and set
$\m{U}$ be equal to $\omega^{\leqslant n}$. $\m{U}$ is naturally
lexicographically ordered. To prove the theorem, we only need to prove
that $\m{U} \arrows{(\m{T})}{\m{S}}{k}$. Indeed, even though $\m{U}$
is not finite, a standard compactness argument can take us to the
finite. 

Let $\{s_i : i < |\m{S}| \}_{<^{\m{S}}_{lex}}$ be a strictly $<^{\m{S}} _{lex}$-increasing enumeration
of the elements of $\m{S}$ and define $f : \funct{|\m{S}|}{|\m{S}|}$ such that
$f(0) = 0$ and $s_{f(i)}$ is the immediate $<^{\m{S}}$-predecessor of
$s_i$ in $\m{S}$ if $i>0$. Define similarly $g : \funct{|\m{T}|}{|\m{T}|}$ for
$\m{T} = \{t_j : j < |\m{T}| \}_{<^{\m{T}}_{lex}}$.
Let also $\mathscr{S} = \{ X \subset \m{U} : X \sqsubset \m{S} \}$
(resp. $\mathscr{T} = \{ X \subset \m{U} : X \sqsubset \m{T}
\}$) where $X \sqsubset \m{S}$ means that $X$ is a $<^{\m{U}}
_{lex}$-initial segment of some $\mc{S} \cong \m{S}$. $\mathscr{S}$ (resp. $\mathscr{T}$) has a natural tree structure with respect to
$<^{\m{U}}
_{lex}$-initial segment, has height $|\m{S}|$ (resp. $|\m{T}|$) and
$\mathscr{S} ^{max}$  (resp. $\mathscr{T}
^{max}$) is equal to $\binom{\m{U}}{\m{S}}$ (resp.
$\binom{\m{U}}{\m{T}}$). Observe that if $X \in \mathscr{S}
\smallsetminus \mathscr{S} ^{max}$ is enumerated as $\{x_i : i < |X| \}_{<^{\m{U}} _{lex}}$ and $u \in \m{U}$ such that $X
<^{\m{U}} _{lex} u$ (that is $x <^{\m{U}} _{lex} u$ for every $x \in X$), then $X \cup
\{ u \} \in \mathscr{S}$ iff $u \in \mathrm{IS}_{\m{U}}(x_{f(|X|)})$ (where
$\mathrm{IS}_{\m{U}}(x)$ denotes the set of immediate
$<^{\m{U}}$-successors of $x$ in $\m{U}$). Consequently, $X, X' \in
\mathscr{S} \smallsetminus \mathscr{S} ^{max}$ can be simultaneously extended
in $\mathscr{S}$ iff $x_{f(|X|)} = x'_{f(|X'|)}$. Now, for $u \in
\m{U}$, let $\mathcal{W} _u$ be a non-principal ultrafilter on
$\mathrm{IS}_{\m{U}}(u)$ and for every $X \in \mathscr{S}
\smallsetminus \mathscr{S}^{max}$, let $\mathcal{V} _X = \mathcal{W}
_{x_{f(|X|)}}$. Hence, $\mathcal{V} _X$ is an ultrafilter on the set
of all elements $u$ in $\m{U}$ which can be used to extend $X$ in
$\mathscr{S}$. Now, let $\mathcal{S}$ be a \textit{$\vec{\mathcal{V}}$-subtree} of
  $\mathscr{S}$, that is a subtree such that for
  every $X \in \mathcal{S} \smallsetminus \mathscr{S} ^{max}$, $\{ u
  \in \m{U} : X <^{\m{U}} _{lex} u \wedge X \cup \{ u \} \in \mathcal{S}
  \} \in \mathcal{V} _X$.    

\begin{claim}
There is $\mc{T} \in \binom{\m{U}}{\m{T}}$ such that
$\binom{\mc{T}}{S} \subset \mathcal{S} ^{max}$.
\end{claim}

For $X \in \mathcal{S}$, let $U_X = \{ u
  \in \m{U} : X <^{\m{U}} _{lex} u \wedge X \cup \{ u \} \in \mathcal{S}
  \}$. $\mc{T}$ is constructed inductively. Start with $\tau _0 =
\emptyset$. Generally, suppose that $\tau _0 <^{\m{U}} _{lex} \ldots
<^{\m{U}} _{lex} \tau _j$ are contructed such that for any $X \subset
\{ \tau _0 , \ldots , \tau _j \}$, if $ X \in \mathscr{S}$, then $X
  \in \mathcal{S}$. Consider now the family $\mathcal{I}$ defined by $\mathcal{I} = \{ I \subset \{ 0, \ldots , j \}
: \{ t_i : i \in I \} \cup \{ t_{j+1} \} \sqsubset \m{S} \}$ and for
$I \in \mathcal{I}$ let $X_I = \{ \tau _i : i \in I \}$. $(X_I)_{I \in
\mathcal{I}}$ is the family of all elements of $\mathscr{S}$ which
need to be extended with $\tau _{j+1}$. In other
words, we have to choose $\tau _{j+1} \in \m{U}$ such that $\{ \tau _0 ,
\ldots , \tau _{j+1} \} \in \mathscr{T}$ and $X_I \cup \{ \tau _{j+1} \} \in
\mathcal{S}$ for every $I \in \mathcal{I}$. To do that, notice that
for any $u \in \m{U}$ which satisfies $\tau _j <^{\m{U}} _{lex} u$, we
have  $\{ \tau _0 ,
\ldots , \tau _j , u \} \in \mathscr{T}$ iff $u \in \mathrm{IS}
_{\m{U}} (\tau _{g(j+1)})$. Now, for any such $u$ and any $I \in
\mathcal{I}$, we have $X_I \cup \{ u \} \in \mathscr{S}$ i.e. $u$
allows a simultaneous extension of all the elements of $\{ X_I : I \in
\mathcal{I} \}$. Consequently, $\mathcal{V} _{X_I}$ does not depend on
$I \in \mathcal{I}$. Let $\mathcal{V}$ be the corresponding common
value. For every $I \in \mathcal{I}$, we have $U_{X_I} \in
\mathcal{V}$ so one can pick $\tau _j <^{\m{U}} _{lex} \tau _{j+1} \in \bigcap_{I \in
  \mathcal{I}} U_{X_I}$. Then $\tau _{j+1}$ is as required. Indeed, on the one
hand  $\{ \tau _0 , \ldots , \tau _{j+1} \} \in \mathscr{T}$ since
$\tau _{j+1} \in \mathrm{IS} _{\m{U}} (\tau _{g(j+1)})$. On the other
hand, for every $I \in \mathcal{I}$, $X_I \cup \{ \tau _{j+1} \} \in
\mathcal{S}$ since $\tau _{j+1} \in U_{X_I}$. At the end of the
construction, we are left with $\mc{T} := \{ \tau _j : j \in |\m{T}|
\} \in \mathscr{T}$ such that $\binom{\mc{T}}{S} \in \mathcal{S}
  ^{max}$, and the claim is proved. 

The proof of the lemma will be complete if we prove the following claim:

\begin{claim}
Given any $k \in \omega \smallsetminus \{ 0 \}$ and
any $\chi : \funct{\binom{\m{U}}{\m{S}}}{k}$, there is a $\vec{\mathcal{V}}$-subtree
$\mathcal{S}$ of $\mathscr{S}$ such that $\mathcal{S} ^{max}$ is
$\chi$-monochromatic.
\end{claim}

We proceed by induction on the height of $\mathscr{S}$. The case
$\mathrm{ht}(\mathscr{S}) = 0 $ is trivial so suppose that the claim
holds for $\mathrm{ht}(\mathscr{S}) = n$ and consider the case
$\mathrm{ht}(\mathscr{S}) = n+1$. Define a coloring $\Lambda :
\funct{\mathscr{S}(n)}{k}$ by $\Lambda (X) = \varepsilon$ iff $\{ u
\in \m{U} : X \cup \{ u \} \in \mathscr{S}(n+1) \wedge \chi (X \cup \{ u
\}) = \varepsilon \} \in \mathcal{V} _X$. By induction hypothesis, we can find a
$\vec{\mathcal{V}}$-subtree $\mathcal{S} _n$ of $\restrict{\mathscr{S}}{n}$ (the tree
formed by the $n$ first levels of $\mathscr{S}$) such that
$\mathcal{S} _n ^{max} $ is $\Lambda$-monochromatic with color
$\varepsilon _0$. This means that for every $X \in \mathcal{S} _n$,
the set 
$V_X := \{ u \in \m{U} : X \cup \{ u \} \in \mathscr{S}(n+1) \wedge \chi (X \cup \{ u
\}) = \varepsilon _0 \}$ is in $ \mathcal{V} _X$. Now, let
$\mathcal{S} = \mathcal{S} _n \cup \{ X \cup \{ u \} : X \in
\mathcal{S} _n \wedge u \in V_X \} $. Then $\mathcal{S}$ is a
$\vec{\mathcal{V}}$-subtree of $\mathscr{S}$ and $\mathcal{S} ^{max}$
is $\chi$-monochromatic.

\end{proof}

We now show how to obtain theorem \ref{thm:2} from lemma
\ref{lemma:2}. Fix $S \subset ]0, + \infty [$, let $(\om{X})$, $(\om{Y}) \in \UU$ and consider $(\ot{T})$ associated to
$(\om{Y})$. As presented in section 2, $(\om{Y})$ can be seen as
$(\textbf{T} ^{max}, <^{\textbf{T}} _{lex})$. Now, notice that there
is a subtree $(\ot{S})$ of $(\ot{T})$ such that for every $(\oc{X})
\in \binom{\textbf{T}^{max} , <^{\textbf{T}} _{lex} }{\om{X}}$, the
downward $<^{\textbf{T}}$-closure of $\mc{X}$ is isomorphic to
$(\ot{S})$. Conversely, for any $(\otc{S}) \in
\binom{\ot{T}}{\ot{S}}$, $(\widetilde{\textbf{S}} ^{max} ,
<^{\widetilde{\textbf{S}}} _{lex})$ is in $\binom{\textbf{T}^{max} ,
  <^{\textbf{T}} _{lex} }{\om{X}}$. These facts allow us to build $(\om{Z})$ such that $ (\om{Z}) \arrows{(\om{Y})}{(\om{X})}{k}$: Apply lemma \ref{lemma:2} to
  get $(\ot{U})$ of height $\mathrm{ht}(\m{T})$ such that
  $(\ot{U}) \arrows{(\ot{T})}{(\ot{S})}{k}$, then simply let $(\om{Z})$ be
  the convexly ordered ultrametric space associated to $(\ot{U})$. To
  check that $(\om{Z})$ works, let $\chi :
  \funct{\binom{\om{Z}}{\om{X}}}{k}$. $\chi$ transfers to $\Lambda :
  \funct{\binom{\ot{U}}{\ot{S}}}{k}$ so find $(\otc{T}) \in
  \binom{\ot{U}}{\ot{T}}$ such that $\binom{\otc{T}}{\ot{S}}$ is
  $\Lambda$-monochromatic. Then the convexly ordered ultrametric space $(\widetilde{\textbf{T}} ^{max}, <^{\widetilde{\textbf{T}}} _{lex})$
 is such that $\binom{\widetilde{\textbf{T}} ^{max},
   <^{\widetilde{\textbf{T}}} _{lex}}{\om{X}}$ is
 $\chi$-monochromatic. But $(\widetilde{\textbf{T}} ^{max},
 <^{\widetilde{\textbf{T}}} _{lex}) \cong (\om{Y})$. Theorem
 \ref{thm:2} is proved.

\section{Ultrametric Urysohn spaces}

Here, $S$ is a countable subset of $]0, + \infty [$. The purpose of this
section is to provide some properties of the Urysohn space $\Ur$
associated to $\U$. $\Ur$ can be seen as follows: The underlying set of $\Ur$
is the set $Q_S$ of all elements $x \in \Q ^S$ which are finitely
supported. As for the distance, it is defined by $d^{\Ur}(x,y) = \min \{s \in S : \forall t \in
S (s<t \rightarrow x(t) = y(t)) \}$. $\Ur$ is really meant to be seen as the set of branches of a
tree. For example, when $S$ is order-isomorphic to $\{1/n : n \in
\omega \}$, $Q_S$ is essentially the set of rational sequences which are
eventually null and $d^{\Ur}$ is the usual distance for the
product topology. With this facts in mind, it is easy to check that $d^{\Ur}$ is an
ultrametric on $\Ur$. Let also $<^{\Ur}_{lex}$ be the natural lexicographical ordering on $\Ur$.

\begin{thm}
$(\Ur , <^{\Ur}_{lex})$ is a countable structures which satisfies i) The finite substructures of $(\Ur , <^{\Ur}_{lex})$ are exactly the
elements of $\UU$ ii) $(\Ur , <^{\Ur}_{lex})$ is ultrahomogeneous, i.e. every isomorphism
between finite substructures of $(\Ur , <^{\Ur}_{lex})$ can be extended to
an automorphism of $(\Ur , <^{\Ur}_{lex})$.
\end{thm}

\begin{proof}
In what follows, we relax the notation and simply write $d$ (resp. $<$)
instead of $d^{\Ur}$ (resp. $<^{\Ur} _{lex}$). i) is easy to check so we concentrate on ii). We proceed by
induction on the size $n$ of the finite substructures.

For $n=1$, if $x$ and $y$ are in $\Ur$, just define $g :
\funct{\Ur}{\Ur}$ by $g(z) = z + y - x$.

For the induction step, assume that the homogeneity of $(\Ur ,
<)$ is proved for finite substructures of size $n$ and
consider two isomorphic substructures of $(\Ur , <)$ of
size $n+1$, namely $x_1 < \ldots < x_{n+1}$ and $y_1 < \ldots <
y_{n+1}$. By induction hypothesis, find $h \in \mathrm{Aut}(\Ur,<)$
such that for every $1 \leqslant i \leqslant n$, $h(x_i) = y_i$. We
now have to take care of $x_{n+1}$ and $y_{n+1}$. Observe first that
thanks to the convexity of $<$, we have $d(x_n,x_{n+1}) = \min \{
d(x_i,x_{n+1}) : 1 \leqslant i \leqslant n \}$ (resp. $d(y_n,y_{n+1}) = \min \{
d(y_i,y_{n+1}) : 1 \leqslant i \leqslant n \}$). Call $s =
d(x_n,x_{n+1}) = d(y_n,y_{n+1})$. Note that $y_{n+1}$ and $h(x_{n+1})$
agree on $S \cap ]s,\infty[$. Indeed, 

\begin{eqnarray*}
d(y_{n+1},h(x_{n+1})) & \leqslant & \max
(d(y_{n+1},y_n),d(y_n,h(x_{n+1}))) \\
& \leqslant & \max (d(y_{n+1},y_n),d(h(x_n),h(x_{n+1})))\\
& \leqslant & \max (s,s) = s
\end{eqnarray*}

Note also that since $y_n < y_{n+1}$ (resp. $h(x_n) < h(x_{n+1})$), we
have $y_n (s) < y_{n+1} (s)$ (resp. $y_n (s) = h(x_n)(s) <
h(x_{n+1})(s)$). So $\R \smallsetminus \Q
\cap ]y_n (s) , \min (y_{n+1}(s), h(x_{n+1})(s)) [$ is non-empty and
has an element $\alpha$. $] \alpha ,
\infty [ \cap \Q$ is order-isomorphic to $\Q$ so we can find a
strictly increasing bijective $\phi : \funct{] \alpha ,
\infty [ \cap \Q}{] \alpha , \infty [ \cap \Q}$ such that \begin{center} $\phi
(h(x_{n+1})(s)) = y_{n+1}(s)$. \end{center} Now, define $j : \funct{\Ur}{\Ur}$ by:

If $d(x,y_{n+1}) > s$ then $j(x) = x$.

If $d(x,y_{n+1}) \leqslant s$ then 

\begin{displaymath}
j(x)(t) = \left \{ \begin{array}{ll}
 x(t) & \textrm{if $t>s$} \\
 x(t) & \textrm{if $t=s$ and $x(t)<\alpha$} \\
 \phi ( x(t)) & \textrm{if $t=s$ and $\alpha < x(t)$}\\
 x(t) + y_{n+1}(t) - h(x_{n+1})(t) & \textrm{if $t<s$} 
 \end{array} \right.
\end{displaymath}

One can check that $j \in \mathrm{Aut}(\Ur ,<)$ and that for every $1
\leqslant i \leqslant n$, $j(y_i) = y_i$. Now, let $g = j \circ h$. We
claim that for every $1 \leqslant i \leqslant n+1$, $g(x_i) =
y_i$. Indeed, if $1 \leqslant i \leqslant n$ then $g(x_i) = j(h(x_i))
= j(y_i) = y_i$. Moreover, 

\begin{displaymath}
g(x_{n+1})(t) = j(h(x_{n+1}))(t) = \left \{ \begin{array}{ll}
 h(x_{n+1})(t) & \textrm{if $t>s$} \\
 \phi (h(x_{n+1})(t)) = y_{n+1}(t) & \textrm{if $t=s$} \\
 h(x_{n+1})(t) + y_{n+1}(t) - h(x_{n+1})(t) = y_{n+1}(t)& \textrm{if $t<s$} 
 \end{array} \right.
\end{displaymath}

i.e. $g(x_{n+1}) = y_{n+1}$.

\end{proof}

We now turn to a description of $\cUr$, the completion of $\Ur$. Note that if $0$ is not an
accumulation point for $S$, then $\Ur$ is discrete and $\cUr =
\Ur$. Hence, in what follows, we will assume that $0$ is an
accumulation point for $S$. 

\begin{thm}
The completion $\cUr$ of the ultrametric space $\Ur$ is the
ultrametric space with underlying set
the set of all elements $x \in \Q ^S$ for which there is
a strictly decreasing sequence $(s_i)_{i \in \omega}$ of elements of
$S$ converging to $0$ such that $x$ is supported by a subset of $\{s_i : i \in \omega \}$.
The distance is given by $d^{\cUr}(x,y) = \min \{s \in S : \forall t \in
S (s<t \rightarrow x(t) = y(t)) \}$.  
\end{thm}

\begin{proof}
We first check that $\Ur$ is dense in $\cUr$. Let $x \in \cUr$ be
associated to the sequence $(s_i)_{i \in \omega}$. For $n \in \omega$,
let $x_n \in \Ur$ be defined by $x_n (s) = x(s)$ if $s>s_n$ and by
$x_n (s) = x(s_n)$ otherwise. Then $d^{\cUr}(x_n , x) = s_{n+1}
\longrightarrow 0$, and
the sequence $(x_n)_{n \in \omega}$ converges to $x$.
To prove that $\cUr$ is complete, let $(x_n)_{n \in \omega}$ be a
Cauchy sequence in $\cUr$. Observe first that given any $s \in S$, the
sequence $x_n (s)$ is eventually constant. Call $x(s)$ the
corresponding constant value.

\begin{claim}
$x \in \cUr$.
\end{claim}

i) is obviously satisfied. To check ii), it is enough to show that
given any $s \in S$, there are $t<s<r \in S$ such that $x$ is null 
on $S \cap ]t,s[$ and on $S \cap ]s,r[$. To do that, fix $t' < s $ in
$S$, and take $N \in \omega$ such that $\forall q \geqslant p
\geqslant N$, $d^{\cUr}(x_q , x_p) < t'$. $x_N$ being in $\cUr$, there
are $t$ and $r$ in $S$ such that $t'<t<s<r$ and $x_N$ is null on $S \cap ]t,s[$ and on $S \cap
]s,r[$. We claim that $x$ agrees with $x_N$ on $S \cap ]t', \infty[$, hence
is null on $S\cap ]t,s[$ and on $S \cap ]s, r[$. Indeed, let $n \geqslant N$. Then $d^{\cUr}(x_n
,x_N) < t' < s$ so $x_n$ and $x_N$ agree on $S \cap ]t', \infty[$. Hence, for every $u \in S \cap
]t',\infty [$, the sequence $(x_n (u))_{n \geqslant N}$ is constant and by
definition of $x$ we have $x(u) = x_n (u)$. The claim is proved. 

\begin{claim}
The sequence $(x_n)_{n \in \omega}$ converges to $x$.
\end{claim}

Let $\varepsilon > 0$. Fix $s \in S \cap ]0 , \varepsilon [$ and $N
\in \omega$ such that $\forall q \geqslant p \geqslant N$, $d^{\cUr}(x_q ,x_p) <
\varepsilon$. Then, as in the previous claim, for every $n \geqslant
N$, $x_n$ and $x_N$ (and hence $x$) agree on $S \cap ]s,
\infty[$. Thus, $d^{\cUr}(x_n , x) \leqslant s < \varepsilon$.

\end{proof}

Since the detailed study of $\cUr$ is not the purpose of this note, we
refer to \cite{B1} for any additional property concerning this space. Let us simply mention that $\cUr$ is
ultrahomogeneous, as well as $(\cUr , <^{\cUr}_{lex})$.

\section{Universal minimal flows}
We now provide some applications of the Ramsey theoretic results proved
in the previous sections to the topological dynamics of isometry
groups of the universal ultrahomogeneous ultrametric spaces presented
in section 5. In this perspective, we start with some general results in topological dynamics
appearing in \cite{KPT}. Let
$G$ be a topological group and $X$ a compact Hausdorff space. A
\textit{$G$-flow} is a continous action $ \funct{G \times X}{X}$. Sometimes,
when is action is understood, the flow is referred to as $X$. Given a
$G$-flow $X$, a nonempty compact $G$-invariant subset $Y \subset X$
defines a subflow by restricting the action to $Y$. $X$ is \textit{minimal}
when $X$ itself is the only nonempty compact $G$-invariant
set (or equivalently, the orbit of any point of $X$ is dense in
$X$). Using Zorn's lemma, it can be shown that every $G$-flow contains
a minimal subflow. Now, given two $G$-flows $X$ and $Y$, a
\textit{homomorphism} from $X$ to $Y$ is a continuous map $\pi : \funct{X}{Y}$
such that for every $x \in X$ and $g \in G$, $\pi ( g \cdot x) = g
\cdot \pi (x)$. An \textit{isomorphism} from $X$ to $Y$ is a bijective
homomorphism from $X$ to $Y$. The following fact is a standard result
in topological dynamics (a proof can be found in \cite{A}):

\begin{thm}    
Let $G$ be a topological group. Then there is a minimal $G$-flow
$M(G)$ such that for any minimal $G$-flow X there is a surjective
homomorphism $\pi : \funct{M(G)}{X}$. Moreover, up to isomorphism,
$M(G)$ is uniquely determined by these properties.  
\end{thm}

$M(G)$ is called the \textit{universal minimal flow} of $G$. Observe that when
$M(G)$ is reduced to a single point, $G$ has a strong fixed point
property: Whenever $G$ acts continuously on a compact Hausdorff space
$X$, there is a point $x \in X$ such that $ g \cdot x = x$ for every
$g \in G$. $G$ is then said to be \textit{extremely amenable}. Now, the results
presented in \cite{KPT} allow to compute the universal minimal flow of
certain groups provided some combinatorial facts hold for a particular
class of finite objects. For example, we saw that $\UU$ is
reasonable (lemma \ref{lemma:1}), has Ramsey property (theorem
\ref{thm:2}) and Ordering property (theorem \ref{thm:3}). There are two corresponding results in topological
dynamics, which read as follows (in the sequel, $\Ur$ is equipped with the discrete topology and
$\cUr$ with the metric topology, whereas transformation groups are
equipped with the corresponding pointwise convergence
topology and the usual composition law):

\begin{thm}
$\mathrm{Aut}(\Ur , <^{\Ur}_{lex})$ is extremely amenable.
\end{thm}

\begin{thm}

\label{thm:8}

The universal minimal flow of $\mathrm{iso}(\Ur)$ is the set $\cLO (\Ur)$ of
convex linear orderings on $\Ur$ together with the action $\funct{\iso (\Ur) \times \cLO (\Ur)}{\cLO (\Ur)}$, $(g,<)
\longmapsto <^g$ defined
by $x <^g y$ iff $g^{-1}(x) < g^{-1}(y)$.
\end{thm}

Let us mention that before \cite{KPT}, there were only very few examples of
non extremely amenable topological groups for which the universal
minimal flows was known to be metrizable, a property that
$M(\iso(\Ur))$ shares. 

In order to obtain analogous results for $\cUr$, we follow the scheme
adopted in \cite{KPT} to prove that the isometry group of
the Urysohn space is extremely amenable. Let $<^{\cUr} _{lex}$ be the natural lexicographical ordering on $\cUr$.

\begin{lemma}

\label{lemma:3}

There is a continuous group morphism under which $\mathrm{Aut}(\Ur , <^{\Ur} _{lex})$ embeds densely into
$\mathrm{Aut}(\cUr , <^{\cUr} _{lex})$.

\end{lemma}

\begin{proof}
Every $g \in \mathrm{iso}(\Ur)$ has unique extension $\hat{g} \in
\mathrm{iso}(\cUr)$. Moreover, observe that $<^{\cUr}_{lex}$ can be reconstituted from
$<^{\Ur} _{lex}$. More precisely, if $\hat{x}, \hat{y} \in \cUr$,
and $x, y \in \Ur$ such that $d^{\cUr}(x,\hat{x}), d^{\cUr}(y,\hat{y}) < d^{\cUr}(\hat{x},\hat{y})$, then
$\hat{x}<^{\cUr} _{lex} \hat{y}$ iff $x<^{\Ur} _{lex}y$. Note that
this is still true when $<^{\cUr} _{lex}$ and $<^{\Ur} _{lex}$ are
replaced by $\prec \in \cLO(\cUr)$ and $\restrict{\prec}{\Ur} \in
\cLO(\Ur)$ respectively. Later, we
will refer to that fact as the \textit{coherence property}. Its first consequence is that the map $g \mapsto \hat{g}$ can actually be seen
as a map from $ \mathrm{Aut}(\Ur , <^{\Ur} _{lex})$ to
$\mathrm{Aut}(\cUr , <^{\cUr} _{lex})$. It is easy to check that it
is a continuous embedding. We now prove that it has dense range. Take $h
\in \mathrm{Aut}(\cUr , <^{\cUr} _{lex})$, $\hat{x}_1 <^{\cUr} _{lex} \ldots
<^{\cUr} _{lex} \hat{x}_n$ in $\cUr$, $\varepsilon > 0$, and consider the
corresponding basic open neighborhood $W$ around $h$. Take $\eta
> 0$ such that $\eta < \varepsilon$ and for every $1 \leqslant i \neq j \leqslant n$, $\eta <
d^{\cUr}(\hat{x}_i , \hat{x}_j)$. Now, pick $x_1 , \ldots , x_n , y_1 ,
\ldots, y_n \in \Ur$ such that for every $1 \leqslant i \leqslant n$, $d^{\cUr}(\hat{x}_i ,
x_i) < \eta$ and $d^{\cUr}(h(\hat{x}_i) , y_i) < \eta$. Then one can
check that the map $x_i \mapsto y_i$ is an isometry
from $\{x_i : 1 \leqslant i \leqslant n \}$ to $\{ y_i : 1 \leqslant i
\leqslant n \}$ (because $\cUr$ is ultrametric) which is also order-preserving
(thanks to the coherence property). By ultrahomogeneity of $(\Ur , <^{\Ur} _{lex})$, we can extend that
map to $g_0 \in \mathrm{Aut}(\Ur , <^{\Ur} _{lex})$. Finally, consider
the basic open neighborhood $V$ around $g_0$ given by $x_1 , \ldots , x_n$
and $\eta$. Then $\{ \hat{g} : g \in V \} \subset W$. Indeed, let $g \in V$. Then: 
\begin{center}
$d^{\cUr}(\hat{g}(\hat{x}_i),h(\hat{x}_i))
\leqslant \max \{ d^{\cUr}(\hat{g}(\hat{x}_i), \hat{g}(x_i)),
d^{\cUr}(\hat{g}(x_i), \hat{g_0} (x_i)), d^{\cUr}(\hat{g}_0 (x_i), h(\hat{x}_i)) \}$
\end{center}

\noindent Now, since $\hat{g}$ is an isometry, $d^{\cUr}(\hat{g}(\hat{x}_i), \hat{g}(x_i)) = d^{\cUr}(\hat{x}_i,
x_i) < \eta < \varepsilon$. Also, since $g \in V$, $d^{\cUr}(\hat{g}(x_i), \hat{g_0} (x_i)) <
\eta < \varepsilon$. Finally, 
$d^{\cUr}(\hat{g}_0 (x_i), h(\hat{x}_i)) = d^{\Ur}(y_i, h(\hat{x_i}))
< \eta < \varepsilon$ by construction of $g_0$. Thus $d^{\cUr}(\hat{g}(\hat{x}_i),h(\hat{x}_i))
< \varepsilon$ and $\hat{g} \in W$.

\end{proof}

\begin{cor}
\label{cor:5}

$\mathrm{Aut}(\cUr ,<^{\cUr} _{lex})$ is extremely amenable.
\end{cor}

\begin{proof}
Let $X$ be a compact Hausdorff space and let $\alpha :
\funct{\mathrm{Aut}(\cUr ,<^{\cUr} _{lex}) \times X}{X}$ be a
continuous action. Then $\beta : \funct{\mathrm{Aut}(\Ur
  ,<^{\Ur} _{lex}) \times X}{X}$ defined by $\beta (g,x) =
\alpha(\hat{g},x)$ is a continuous action. Since
$\mathrm{Aut}(\cUr ,<^{\cUr} _{lex})$ is extremely amenable, there is
$x \in X$ which is fixed under $\beta$. Now, $\{ \hat{g} : g \in
\mathrm{Aut}(\Ur ,<^{\Ur} _{lex})\}$ being dense in $\mathrm{Aut}(\cUr
,<^{\cUr} _{lex})$, $x$ is also fixed under $\alpha$.

\end{proof}

Let us now look at the topological dynamics of the isometry group
$\iso (\cUr)$. Note that $\iso (\cUr)$ is not extremely amenable as
its acts continuously on the space of all convex linear orderings
$\cLO (\cUr)$ on $\cUr$ with no fixed point. The following result
shows that in fact, this is its universal minimal compact action.

\begin{cor}
The universal minimal flow of $\mathrm{iso}(\cUr)$ is the set $\cLO
(\cUr)$ of convex linear orderings on $\cUr$ together with the action $\funct{\iso (\cUr) \times \cLO (\cUr)}{\cLO (\cUr)}$, $(g,<)
\longmapsto <^g$ defined
by $x <^g y$ iff $g^{-1}(x) < g^{-1}(y)$.
\end{cor}

\begin{proof}

Equipped with the topology for which the basic open sets
are those of the form $\{ \prec \in \cLO (\cUr) : \restrict{\prec}{X}
= \restrict{<}{X} \}$ (resp. $\{ \prec \in \cLO (\Ur) : \restrict{\prec}{X}
= \restrict{<}{X} \}$) where $X$ is a finite subset of $\cUr$
(resp. $\Ur$), the space $\cLO (\cUr)$ (resp. $\cLO (\Ur)$) is compact. To see that the action is continuous,
let $< \in \cLO (\cUr)$, $g \in \iso (\cUr)$ and $W$ a basic open
neighborhood around $< ^g$ given by a finite $X \subset \cUr$. Now take $\varepsilon >
0$ strictly smaller than any distance in $X$
and consider $U = \{ h \in \iso (\cUr) : \forall x \in X
(d^{\cUr}(g^{-1}(x),h^{-1}(x)) < \varepsilon) \}$. Let also $V = \{ \prec \in  \cLO (\cUr) :
\restrict{\prec}{\overleftarrow{g}X} =
\restrict{\prec}{\overleftarrow{h}X} \}$ where $\overleftarrow{g}X$
(resp. $\overleftarrow{h}X$) 
denotes the inverse image of $X$ under $g$ (resp. h). We claim that
for every $(h, \prec) \in U \times V$, we have $\prec ^h \in W$. To
see that, observe first that if
$x,y \in X$, then $h^{-1}(x) \prec
h^{-1}(y)$ iff $g^{-1}(x) \prec g^{-1}(y)$ (this is a consequence of
the coherence property). So if $(h, \prec) \in U \times V$ and $x,y \in X$ we have  

\begin{eqnarray*}
x \prec ^h y & \textrm{iff} & h^{-1}(x) \prec h^{-1}(y) \quad
\textrm{by definition of} \, \,\prec ^h \\
& \textrm{iff} & g^{-1}(x) \prec g^{-1}(y) \quad \textrm{by the
  observation above} \\
& \textrm{iff} & g^{-1}(x) < g^{-1}(y) \quad
\textrm{since} \, \, h \in U \\
& \textrm{iff} & x <^g y \quad \textrm{by definition of} <^g
\end{eqnarray*}

So $\prec ^h \in W$ and the action is continuous. To prove the theorem, we have to show 

\begin{claimm}

\label{claimm:1}

Given any $< \in
\cLO(\cUr)$, the $\iso (\cUr)$-orbit of $<$ is dense in $\cLO
(\cUr)$.
\end{claimm}

\begin{claimm}

\label{claimm:2}
Given any minimal $\iso (\cUr)$-flow $X$, there is a continuous and onto $\alpha : \funct{\cLO (\Ur)}{X}$ such that: 
$\forall g \in \iso (\Ur)$, $\forall x \in X$, $\alpha (\varphi (g)
\cdot x) = \varphi (g) \cdot \alpha (x)$.
\end{claimm}

For claim \ref{claimm:1}, let $\prec \in \cLO (\cUr)$ and $X$ be a finite subset of
$\cUr$. We want to produce $g \in \iso (\cUr)$ such that $<^g$ and
$\prec$ agree on $X$. To do that, let $\psi$ be the restriction map
defined by $\psi : \funct{\cLO (\cUr)}{\cLO (\Ur)}$
with $\psi (<) = \restrict{<}{\Ur}$. Now, take $\varepsilon > 0$ strictly
smaller than all the distances in $X$ and for every $x \in X$, let $x'
\in \Ur$ be such that $d(x,x') < \varepsilon$ (here and later $d$ will stand
for $d^{\cUr}$). Thanks to the coherence property,
$\forall x, y \in X$, $x < y$ (resp. $x \prec y$) iff $x' < y'$
(resp. $x' \prec y'$). Call $X' = \{x' : x \in X \}$. Now, $\psi (<)$ and $\psi (\prec)$ are
in $\cLO (\Ur)$; since $\cLO(\Ur)$ is a minimal $\iso (\Ur)$-flow, the
$\iso (\Ur)$-orbit of $\psi (<)$ is dense in $\cLO(\Ur)$ so there is $h
\in \iso (\Ur)$ such that $\psi(<)^h$ and $\psi(\prec)$ agree on
$X'$. Then $g = \hat{h}$ is such that $<^g$ and $\prec$ agree on
$X$: Let $x,y \in X$. Then $g(x) < g(y)$ iff $g(x') < g(y')$ iff
$h(x') \, \, \psi (<) \, \, h(y')$ iff $ x' \, \, \psi( \prec) \, \,
y' $ iff $x' \prec y'$ iff $x \prec y$ and claim \ref{claimm:1} is proved.

For claim \ref{claimm:2}, consider a minimal $\iso (\cUr)$-space $X$
and observe first that there is a natural dense embedding $\varphi : \funct{\iso
(\Ur)}{\iso (\cUr)}$ (recall that $\iso (\Ur)$ is equipped with the
pointwise convergence topology coming from the discrete topology on
$\Ur$ whereas $\iso (\cUr)$ is equipped with the pointwise
convergence topology coming from the metric topology on
$\cUr$). $\varphi$ allows to consider the action $\funct{\iso
  (\Ur) \times X}{X}$ defined by $\beta (g,x) = \varphi (g) \cdot
x$. One can check that this action is continuous so since $ \cLO
(\Ur)$ is the universal minimal flow of $\iso (\Ur)$,
there is $\pi : \funct{\cLO (\Ur)}{X}$ continuous and onto such that: 
$\forall g \in \iso (\Ur)$, $\forall x \in X$, $\pi (\varphi (g)
\cdot x) = \varphi (g) \cdot \pi (x)$. Now, let $\widehat{\pi} :
\funct{\cLO (\cUr) }{X}$ be defined as $\pi \circ \psi$. We claim that
$\widehat{\pi}$ is continuous and onto. To see that, it is enough to
show that $\psi$ is continuous and onto. Surjectivity comes from the fact that any element $<$ of
$\cLO (\Ur)$ extends naturally to $\widehat{<} \in \cLO (\cUr)$: $x
\widehat{<} y$ iff there are $x',y' \in \Ur$ such that $d^{\cUr}(x,x'),
d^{\cUr}(y,y') < d^{\cUr}(x,y)$ and $x' < y'$. For
continuity, let $< \in \cLO (\Ur)$ and $X$ a finite subset of
$\Ur$. Consider the basic open neighborhood $V$ around $<$ in
$\cLO(\Ur)$ defined by $V = \{ < ' \in \cLO (\Ur)
: \restrict{< '}{X} = \restrict{<}{X} \}$. Then every element of $\{ \prec \in \cLO (\Ur)
: \restrict{\prec}{X} = \restrict{<}{X} \}$ is such that $\psi (\prec)
\in V$. So to show that $\alpha = \widehat{\pi}$ works and complete the proof, it remains to prove that 

\begin{center}
$\forall g \in \iso (\cUr) \, \forall < \in \cLO (\cUr) \quad
\widehat{\pi} (<^g) = \widehat{\pi}(<) ^g$
\end{center}

Let $g \in \iso (\cUr)$ and $< \in \cLO (\cUr)$. Writing $D$ for
$\mathrm{ran}(\varphi)$, observe first that thanks to the continuity
of the action $\funct{\iso (\cUr) \times \cLO (\cUr)}{\cLO (\cUr)}$ and
the fact that $D$ is dense, we have  
\[<^g \, = \lim _{ \substack{h \rightarrow g \\ h \in D}} <^h \]

Thus, \[ \widehat{\pi}(<^g) = \widehat{\pi}(\lim _{ \substack{h \rightarrow g \\ h \in D }}
 <^h) = \lim _{ \substack{h \rightarrow g \\ h \in D}}
 \widehat{\pi}(<^h) = \lim _{ \substack{h \rightarrow g \\ h \in D }}
 \pi (\psi(<^h)) \]

Now, notice that $\psi (<^h) = \psi (<) ^h$. Indeed,
let $x, y \in \Ur$. Then 
\begin{eqnarray*}
x \, \, \psi (<^h) \, \, y & \textrm{iff} & x <^h y \quad \textrm{since} \, \, x, y \, \in \Ur\\
& \textrm{iff} & h^{-1}(x) < h^{-1}(y) \quad \textrm{by definition of} \, \,<^h \\
& \textrm{iff} & h^{-1}(x) \, \, \psi (<) \, \, h^{-1}(y) \quad \textrm{since} \,
\, h^{-1}(x), h^{-1}(y) \, \in \Ur \\
& \textrm{iff} & x \, \, \psi (<) ^h \, \, y \quad
\textrm{by definition of} \, \, \psi (<) ^h 
\end{eqnarray*}

So \[ \widehat{\pi}(<^g) = \lim _{ \substack{h \rightarrow g \\ h \in D }}
 \pi (\psi(<^h)) = \lim _{ \substack{h \rightarrow g \\ h \in D }}
 \pi (\psi(<) ^h) = \lim _{ \substack{h \rightarrow g \\ h \in D }}
 \pi (\psi(<)) ^h = \widehat{\pi} (<) ^g\]
 
\end{proof}

We finish with two direct consequences of the previous corollary. The
first one is a purely topological comment along the lines of the
remark following theorem \ref{thm:8}: To show that the underlying
space related to the universal minimal flow of $\iso(\cUr)$ is $\cLO (\cUr)$, we use the fact
that the restriction map $\psi : \funct{\cLO (\cUr)}{\cLO (\Ur)}$
defined by $\psi (<) = \restrict{<}{\Ur}$ is continuous and
surjective. However, it turns out that $\psi$ is actually a
homeomorphism. $\cLO (\Ur)$ being metrizable, we consequently get:

\begin{cor}
The underlying space of the universal minimal flow of $\iso(\cUr)$ is metrizable.
\end{cor}

The second consequence is based on the simple observation that when
the distance set $S$ is $\{1/n : n \in \omega \}$, $\cUr$ is the Baire
space $\mathcal{N}$. Hence:

\begin{cor}
When $\mathcal{N}$ is equipped with the product metric, the universal
minimal flow of $\mathrm{iso}(\mathcal{N})$ is the set of all convex
linear orderings on $\mathcal{N}$. 
\end{cor}

\section{Remarks for further studies}
As written in the introduction, very little is known about Ramsey
properties for classes of finite metric spaces so there is a lot to do
in this direction. Unfortunately, we doubt that the generalization
can be pushed much further. Indeed, the structural connection with
trees is probably too specific to be representative of the generic
behaviour of finite metric spaces. For example, if the extreme
amenability of the unitary group of $\ell _2$ can be proved via the
approach of \cite{KPT}, there is little hope that it can be done with
the techniques of this article. Finally, in section 4,
we presented a Ramsey result concerning a class
of finite ordered trees with fixed finite height. Our hope is that the
method we adopted may be used for the partition calculus of countable
trees.

\end{document}